# Analytic exploration of safe basins in a benchmark problem of forced escape

*Gleb Karmi, Pavel Kravetc and Oleg Gendelman*

Faculty of Mechanical Engineering, Technion – Israel Institute of Technology, Haifa, 3200003, Israel

## Abstract

The paper presents an analytical approach for predicting the *safe basins* (SB) in a plane of *initial conditions* (IC) for escape of classical particle from the potential well under harmonic forcing. The solution is based on the approximation of isolated resonance, which reduces the dynamics to conservative flow on a two-dimensional *resonance manifold* (RM). Such a reduction allows easy distinction between escaping and non-escaping ICs. As a benchmark potential, we choose a common parabolic-quartic well with truncation at varying energy levels. The method allows accurate predictions of the SB boundaries for relatively low forcing amplitudes. The derived SBs demonstrate an unexpected set of properties, including decomposition into two disjoint zones in the IC plane for a certain range of parameters. The latter peculiarity stems from two qualitatively different escape mechanisms on the RM. For higher forcing values, the accuracy of the analytic predictions decreases to some extent due to the inaccuracies of the basic isolated resonance approximation, but mainly due to the erosion of the SB boundaries caused by the secondary resonances. Nevertheless, even in this case the analytic approximation can serve as a viable initial guess for subsequent numeric estimation of the SB boundaries.

## 1. Introduction and historical survey

Escape from a potential well is a tool for modeling a plethora of effects in many branches of science and engineering [1]–[6]. Examples of the escape phenomena include gravitational collapse in celestial mechanics [7], energy harvesting [6], absorption of particles [2], [8], responses of Josephson junctions [8], resonance capture [9], [10], capsize of ships [3], [11]–[13], dynamic pull-in in microelectromechanical systems [14], collapses of arches [4], energy transport [15], and capture of charged particles by an electric field [16] to mention a few.

One of the first seminal works on escape was written by Kramers almost 80 years ago [17]. While focusing on thermal activation of chemical reactions, it created a multitude of novel research areas. A particularly interesting phenomenon is the escape triggered by a stochastic resonance [18]–[20].

In recent years a lot of attention was drawn to the dynamics of escape under the external harmonic forcing [21]–[26]. In [21] the escape from a squared hyperbolic secant potential was treated with the help of approximation of isolated 1:1 *resonance* (AIR). The AIR technique consists of two main steps: the *action-angle* (AA) transformation and the subsequent averaging over the fast phases. As a result, for the Hamiltonian system it yields an integral of motion describing a family of *resonance manifolds* (RM), conveniently described by a two-dimensional phase portrait of the slow flow. Qualitative analysis of the slow-flow equation allows to identify the escaping trajectories, to describe the mechanisms of transition to escape and to find the critical values of parameters.

Previous studies considered only special (primarily, zero) ICs, therefore, the question about the global dynamics remains open. This question can be formulated as follows: what is the set of all ICs corresponding to non-escape for a potential well? The aforementioned set of the ICs for the non-escaping trajectories is commonly referred to as a *safe basin* (SB) [27], [28]. In order to quantify the size of an SB the notion of integrity measure (e.g., *global integrity measure* (GIM) which is the hyper-volume of the SB) is introduced [27]–[29]. Reduction of the integrity with the change of the system parameters is referred to as erosion, and the variation of integrity is called the erosion profile. In [27] the control was applied by alteration of the excitation shape leading to the shift of the erosion profile towards greater external forcing amplitudes by eliminating homo/heteroclinic bifurcations of hilltop saddles. Soliman and Thompson investigated *stochastic integrity* (SI) in the presence of a white noise [30] while Orlando et al. continued the study SI for bistable elastic structures with harmonic and stochastic loading [31], [32].

Despite all the advances in the study of SBs and their erosion, most of the related methods have significant limitations. First, due to the complexity of the problem, the analytic treatment of SBs often stays beyond the consideration. Secondly, an integrity measure only quantifies the magnitude of an SB, however, it possesses no information about the shape. The boundary for possible placement of the attractors in a strongly nonlinear forced system can be assessed by semi-analytic procedure [33] and can be highly nontrivial.

Current work suggests the analytic approach for predicting the SB boundaries in a benchmark model system — a particle in truncated quartic potential — under harmonic forcing with constant amplitude. First, the AIR is applied to find the escape threshold for arbitrary initial conditions inside the well, under different levels of truncation. Then, the SB boundaries are identified and classified both on the RM and in the IC plane.

The paper is organized as follows. In Section 2 the model is presented, and the AIR is described. Those are followed by derivation of the critical forcing amplitude and classification of the safe basins. Section 3 presents the comparison of analytic results with numerical simulations and addresses the limitations of the suggested approach. Section 5 is devoted to discussion and conclusions.

## 2. Model description and analytic results

We consider a harmonically forced quartic oscillator in the absence of damping. In the nondimensional form the equation is written as

$$\ddot{q}(\tau) + q(\tau) - q(\tau)^3 = F\sin(\tau\Omega + \Psi), \tag{1}$$

where $q$ is the displacement, $\tau$ stands for a nondimensional time, $F$ is the amplitude of the external force, $\Omega$ is the force frequency, and $\Psi$ is the forcing phase.

Equation (1) defines motion of a harmonically forced particle of the unit mass in a symmetric quartic potential

$$V(q) = \frac{q^2}{2} - \frac{q^4}{4}. \tag{2}$$

Potential (2) has a well (local minimum) with barriers (local maxima) at $q = \pm 1$. Assuming initially the particle is trapped inside the potential well, we are looking for the trajectories that correspond to the escape, i.e., particle leaving the well. Commonly, the escape is defined as $\lim_{\tau \to \infty} |q(\tau)| > q_{\max}$ where $0 < q_{\max} \leq 1$. However, this criterion is problematic for numerical implementation, therefore, we use an alternative definition of escape: $\max_{\tau} |q(\tau)| > q_{\max}$. Note that the two definitions are not equivalent, nonetheless, for the purpose of our problem, the latter definition is sufficient [21], [23].

In light of the range constraints due to the escape, one can generalize the potential (2) by applying the truncation at certain critical displacement inside the potential well:

$$U(q) = \begin{cases} \dfrac{q^2}{2} - \dfrac{q^4}{4} - \dfrac{q_{\max}^2}{2} + \dfrac{q_{\max}^4}{4}, & |q| < q_{\max} \\ 0, & |q| \geq q_{\max} \end{cases}. \tag{3}$$

Alternatively, we can perform the truncation by setting the maximum energy level $\xi_{\max}$. In this case the maximum displacement $q_{\max}$ is defined as follows:

$$q_{\max} = \sqrt{1 - \sqrt{1 - 4\xi_{\max}}}. \tag{4}$$

At a first glance, the truncated potential may appear far-fetched, however, in some applications it emerges naturally. For example, in the context of naval mechanics the escape from the full (non-truncated) potential corresponds to the ship capsizing [18], while the truncation corresponds to additional constraints imposed on the maximum allowed roll angle of the ship. Furthermore, truncated polynomial potentials can be used as a tool for approximating the escape problem in more intricate potentials.

It is convenient to represent equation (1) in the Hamiltonian form

$$\dot{q} = -\frac{\partial H}{\partial p}, \quad \dot{p} = \frac{\partial H}{\partial q}, \tag{5}$$

where $p = \dot{q}$ is the momentum of the particle and Hamiltonian $H(p,q)$ is

$$H(p,q) = H_0(p,q) - Fq\sin(\Omega\tau + \Psi). \tag{6}$$

Here, $H_0$ denotes the basic Hamiltonian

$$H_0(p,q) = \frac{p^2}{2} - \frac{q^4}{4} + \frac{q^2}{2}, \quad -q_{\max} < q < q_{\max}, \tag{7}$$

which describes the free motion of the particle. The basic Hamiltonian $H_0$ provides us with yet another definition of escape: $\max_\tau H_0(q(\tau), p(\tau)) > \xi_{max}$ where $0 < \xi_{max} \leq 1/4$, i.e., the complete energy exceeds a maximum level $\xi_{max}$. The difference between different escape definitions is discussed in Appendix.

## 2.1 Action angle formalism and isolated resonance approximation

In order to obtain isolated resonance approximation, one needs to rewrite the Hamiltonian system (4) in AA variables and subsequently average it over the fast phases.

The canonical pair $(I, \theta)$ composed of action and angle is defined as follows [1]:

$$I = \frac{1}{2\pi} \oint_{\Gamma_E} p(q, E) dq, \quad \theta = \frac{\partial}{\partial I} \int_0^q p(x, I) dx, \qquad (8)$$

where $\Gamma_E$ is a curve defined by a level set $H_0 = E$. The canonical transformation (8) does not explicitly depends on time $\tau$, and therefore, the Hamiltonian (6) can be written in AA variables:

$$H(I, \theta) = H_0(I) - Fq(I, \theta) \sin(\Omega \tau + \Psi). \qquad (9)$$

Due to the $2\pi$-periodicity of the angle variable, the Hamiltonian (9) can be expanded in terms of the Fourier series:

$$H(I, \theta) = H_0(I) - \frac{iF}{2} \sum_{m=-\infty}^{\infty} q_m(I)(e^{i(m\theta + \Omega\tau + \Psi)} - e^{i(m\theta - \Omega\tau - \Psi)}), \quad q_m = q_{-m}^*, \qquad (10)$$

where symbol $*$ denotes the complex conjugation. In order to treat 1:1 resonance the phase $\gamma = \theta - \Psi - \tau\Omega$ is assumed to be slow, while all other combinations of phases must be considered fast [15], [21], [23]. After performing the averaging over the fast phases, one arrives at the following slow-flow equations:

$$\dot{J} = -\frac{\partial C(\gamma, J)}{\partial \gamma}, \quad \dot{\gamma} = \frac{\partial C(\gamma, J)}{\partial J}, \qquad (11)$$

with the first integral [15], [21]

$$C(\gamma, J) = H_0(J) - \frac{iF}{2}(q_1(J)e^{i\gamma} - q_1^*(J)e^{-i\gamma}) - \Omega J = const, \qquad (12)$$

where $J = \langle I(\tau) \rangle$ denotes the averaged action. Equation (12) defines a family of 1:1 RMs of the system. The existence of the first integral (12) is attributed to the fact that the system (5) is Hamiltonian. The value of the constant is determined by the IC on the RM, i.e., the values of the averaged action $J$ and the slow phase $\gamma$ at which the system is captured by the RM.

By applying transformation (8) to the basic Hamiltonian (7) one can express action-angle variables as functions of the energy $E$ and the displacement $q$ as follows (see [23] for the detailed derivation):

$$I(E) = \frac{2\sqrt{2}}{3\pi}\sqrt{1+\mu(E)}(\mathbf{E}(k(E)) - \mu(E)\mathbf{K}(k(E))), \tag{13}$$

$$\theta(q,E) = \frac{\pi}{2\mathbf{K}(k(E))}\mathbf{F}\left(\sin^{-1}\left(\frac{q}{\sqrt{1-\mu(E)}}\right)\bigg|k(E)\right), \tag{14}$$

where $\mu(E) = \sqrt{1-4E}$; functions $\mathbf{K}(k)$ and $\mathbf{E}(k)$ are complete elliptic integrals of the first and the second kind, respectively, with the modulus[1] $k = \sqrt{(1-\mu)/(1+\mu)}$. Finally, $\mathbf{F}(\varphi, k)$ denotes the incomplete elliptic integral of the first kind.

Equation (14) yields the displacement $q$ as a function of the angle $\theta$ and the energy $E$:

$$q(\theta, E) = \sqrt{1-\mu(E)}\,\text{sn}\left(\frac{2\theta K(k(E))}{\pi}\bigg|k(E)\right), \tag{15}$$

where sn($x, k$) denotes the Jacobi elliptic sine function.

By taking the time derivative of (15) and noting that

$$\dot{\theta}(E) = \frac{dE}{dI} = \frac{\pi\sqrt{1+\mu(E)}}{2\sqrt{2}\mathbf{K}(k(E))}, \tag{16}$$

we obtain the following expression for the momentum $p$:

$$p(\theta, E) = \sqrt{\frac{1-\mu(E)^2}{2}}\,\text{cn}\left(\frac{2\theta\mathbf{K}(k(E))}{\pi}\bigg|k(E)\right)\text{dn}\left(\frac{2\theta\mathbf{K}(k(E))}{\pi}\bigg|k(E)\right). \tag{17}$$

In order to proceed with derivation of the conservation law (12), one needs to find $H_0(I)$. Unfortunately, the inversion of (13) in a closed form is impossible, however, for the purpose of the escape problem it is enough to consider an averaged energy $\xi = \langle E(t)\rangle$ instead of the averaged action $J$. In this case the conservation law for the slow dynamics is written as follows:

---

[1] Commonly in scientific software like SciPy or Wolfram Mathematica one needs to square the modulus, i.e., $\mathbf{K}(k)$ or $\mathbf{F}(\phi|k)$ in the formula reads as $\mathbf{K}(k^2)$ or $\mathbf{F}(\phi|k^2)$ in the program code.

$$C(\xi,\gamma) = -FG(\xi)\cos(\gamma) - \Omega J(\xi) + \xi = const, \tag{18}$$

where

$$J(\xi) = \frac{2\sqrt{2}}{3\pi}\sqrt{1+\mu(\xi)}\,(\mathbf{E}(k(\xi)) - \mu(\xi)\mathbf{K}(k(\xi))), \tag{19}$$

$$G(\xi) = \frac{\pi\sqrt{1+\mu(\xi)}}{\mathbf{K}(k(\xi))}\exp\left(-\frac{\pi\mathbf{K}\left(\sqrt{1-k(\xi)^2}\right)}{2\mathbf{K}(k(\xi))}\right). \tag{20}$$

## 2.2 Escape threshold on the control plane

We proceed with the qualitative analysis of the averaged dynamics of the system defined by the integral of motion (18). In order to determine if a trajectory starting from a given IC ($q_{ini}, p_{ini}$) is escaping according to the maximum energy criterion, one has to check whether the corresponding level curve

$$C(\gamma,\xi) = C(\gamma_{ini},\xi_{ini}), \tag{21}$$

crosses the circle $\xi = \xi_{max}$. In [21], [23] two main scenarios of transition to escape were identified as *a maximum mechanism* (MM) and *saddle mechanism* (SM).

The maximum mechanism is a passage to escape described as follows. For a fixed frequency $\Omega$, and the forcing $F < F_{cr}$ the phase curve (21) lies entirely below the threshold level $\xi = \xi_{max}$ and the particle stays in the well; at $F = F_{cr}$ the curve (21) is tangent to the circle $\xi = \xi_{max}$; finally, for values $F > F_{cr}$ the same curve crosses the circle $\xi = \xi_{max}$, leading to the escape.

Similarly, the saddle mechanism represents the following behaviour. For the values of the external force $F < F_{cr}$, the part of the curve (21) connected to the point ($\gamma_{ini}, \xi_{ini}$) remains below the saddle on the RM; when $F = F_{cr}$ it passes through the saddle point; and for $F > F_{cr}$ the escape occurs.

In [23] both mechanisms are described and illustrated for the quartic potential (2) and the limiting phase trajectory (LPT), i.e., $C(\gamma, \xi) = 0$.

In order to obtain the critical value $F_{cr}$ of the external forcing amplitude needed for the escape through MM, one can solve the following equation:

$$C(\gamma^*, \xi_{max}) = C(\gamma_{ini}, \xi_{ini}). \tag{22}$$

where $\gamma^*$ is such that

$$\left.\frac{\partial C}{\partial \gamma}\right|_{\gamma = \gamma^*, \xi = \xi_{max}} = 0. \qquad (23)$$

It is easy to see that for any $0 < \xi_{max} < 1/4$ the possible values of $\gamma^*$ are 0 and $\pi$. Therefore, the critical force $F_{cr}$ is

$$F_{cr} = \pm \frac{C_0 + \Omega J(\xi_{max}) - \xi_{max}}{G(\xi_{max})}. \qquad (24)$$

In case of SM, the critical force $F_{cr}$ is obtained by solving the system

$$\nabla C\left(\gamma^\dagger, \xi^\dagger\right) = 0, \quad \det H(C) < 0, \qquad (25)$$

$$C\left(\gamma^\dagger, \xi^\dagger\right) = C(\gamma_{ini}, \xi_{ini}), \qquad (26)$$

where $\nabla$ and $H$ denote the gradient and the Hessian matrix, respectively. Note that, solving the above equations analytically for $\xi^\dagger$ is an impossible task. However, it is not necessary as one can use $\xi^\dagger$ to parameterize the curve $\left(\Omega\left(\xi^\dagger\right), F\left(\xi^\dagger\right)\right)$. Figures 1 and 2 show the behaviour of the critical force curve depending on different values of $\xi_{ini}$ and $\gamma_{ini}$, respectively.

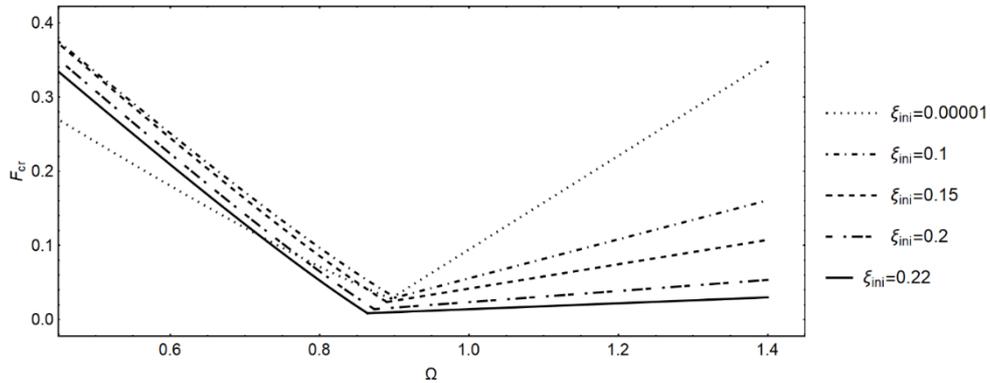

**Figure 1**: Escape threshold versus excitation frequency for various values of $\xi_{ini}$ and $\gamma_{ini}$=0.25, $\xi_{max}$=0.242

Generally, higher initial energies $\xi_{ini}$ correspond to lower critical force amplitude required for the escape. However, our research shows that zero initial energy in the case of SM turns out to be an exception as its $F_{cr}(\Omega)$ curve lies below analogous curves for some higher initial energies.

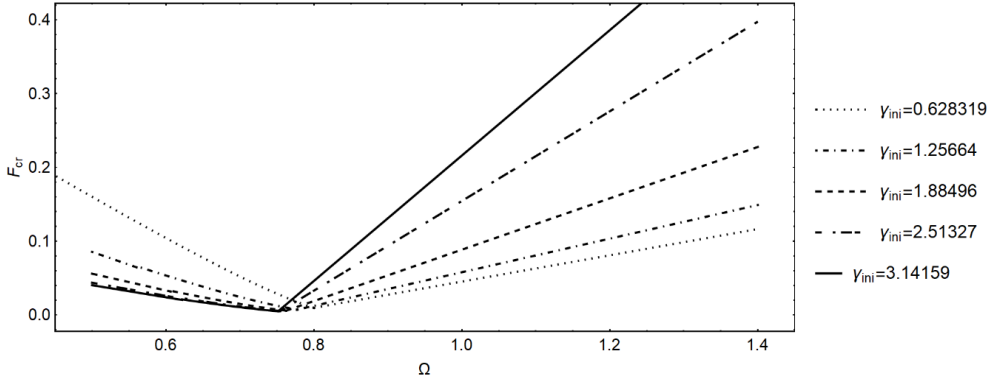

**Figure 2**: Escape threshold versus excitation frequency for various values of $\gamma_{ini}$ and $\xi_{ini}=0.15$, $\xi_{max}=0.242$

Figure 2 demonstrates critical force curves $F_{cr}(\Omega)$ for various values of initial angle $\gamma_{ini}$. As one can easily see that increasing the value of $\gamma_{ini}$ lowers the $F_{cr}$ curve for the SM, while raising the corresponding curve for the MM. Note that, the slow phase $\gamma$ contains the phase $\Psi$ of the external forcing. Its effect on the global dynamics is demonstrated in Section 3.2.

Figure 3 and Figure 4 show curves $F_{cr}(\Omega)$ for the zero IC and various values of the energy truncation level $\xi_{max}$. In particular, Figure 3 demonstrates three different escaping mechanisms: green dashed line corresponds to the MM through the tangent point ($\gamma = \pi, \xi = 0.15$), blue solid curve corresponds to the SM, and red dotted line — to the MM through the point ($\gamma = 0, \xi = 0.15$). The latter scenario occurs when the energy $\xi^{\dagger}$ at the saddle point is above the truncation level $\xi_{max} = 0.15$, and therefore, we will reasonably call this particular case of the MM a *saddle maximum mechanism* (SMM).

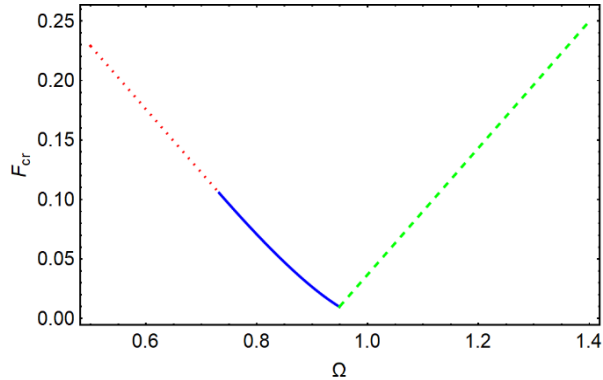

**Figure 3**: Escape threshold versus excitation frequency for the zero IC and the truncation energy level $\xi_{max} = 0.15$. Dotted, solid and dashed lines correspond to SMM, SM and MM, respectively.

Figure 4 presents the transformation of the critical forcing curve $F_{cr}(\Omega)$ with increasing depth of the well. The shallower the well, the lower the forcing amplitude is, while the resonance frequency (at the minima) increases. This behaviour is attributed to the fact that a shallow well can be approximated with a quadratic potential (i.e., linear oscillator), hence the resonance frequency is close to $\Omega = 1$.

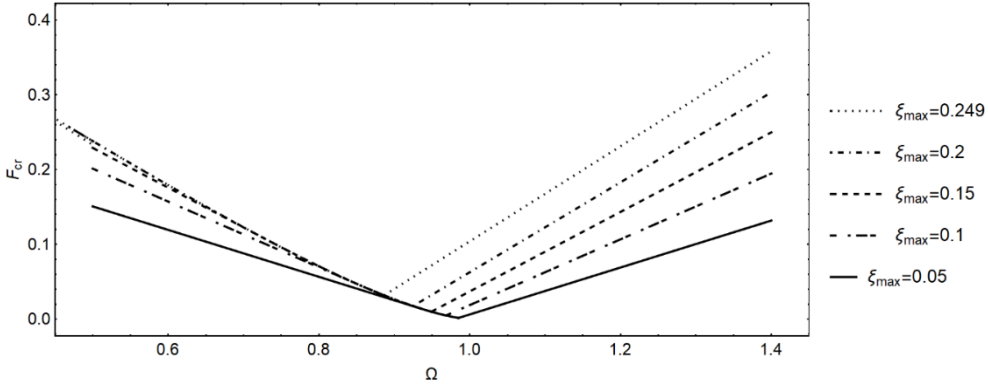

**Figure 4**: Escape threshold versus excitation frequency for the zero IC for various values of the threshold energy level $\xi_{max}$

## 2.3 Safe basins in the vicinity of the main resonance of the perturbed system

The escape of the particle corresponds to the phase curves that transverse the threshold energy level $\xi = \xi_{max}$ on the phase cylinder $(\gamma, \xi)$. On the contrary, the phase curves located below the circle $\xi = \xi_{max}$ correspond to non-escaping trajectories of the particle. Therefore, in order to locate the SB, one needs to determine its boundary — the critical position of the non-escaping phase curves.

Similarly, to the available escape mechanisms, there are two types of the SB depending on their boundaries: the maximum type and the saddle type. We say the basin is of *maximum type* (SBMT) if its boundary is a phase curve tangent to the circle $\xi = \xi_{max}$. An SB is of the *saddle type* (SBST) if its boundary passes through the saddle point. Note that these two types are not mutually exclusive and can coexist.

The boundary of a SBMT is defined by the equation:

$$C(\gamma,\xi)=C(\gamma^*,\xi^*), \qquad (27)$$

where $\gamma^*$ is a solution to (23). It is clear that there are only two possible values of $\gamma^*$: $\gamma^* = 0$ and $\gamma^* = \pi$. Once the boundary is established, one needs to locate the safe region it is enclosing. There are two topologically non-equivalent kinds of SBMT: "island" and "peninsula". The main difference between the two is that the peninsula wraps around the RM phase cylinder, whereas the island does not. However, for our purposes the most important distinction is in the way we find the safe region. In the case of an island the safe region is defined as all the points enclosed by the curve (28) whereas in the case of peninsula the safe region comprises all the points bounded between the curve (28) and the circle $\xi=0$. It is important to note that the latter kind of SBMT contains points with all possible values of the angle $\gamma$. A special case of an island SBMT has an LPT as its boundary. In this case the safe region is defined the same way as for a peninsula.

It is important to remark that there is a singularity at $\xi=1/4$, thus the tangent point $(\gamma^*, \xi^*=1/4)$ is not defined.

After the boundary of SB in the space $(\gamma, \xi)$ is computed, one can use formulas (15), (17) to obtain the corresponding SB in the space of IC on the $(q_0,p_0)$ plane.

The island kind of SBMT is illustrated in Figure 5, while the examples of the peninsula SBMT are shown in Figure 6 and Figure 7. The left panels show the level curves of the first integral (18), maximum energy level (blue dashed line) and the SB (grey region with red boundary). Panels on the right present the SB on the $(q_0,p_0)$ plane.

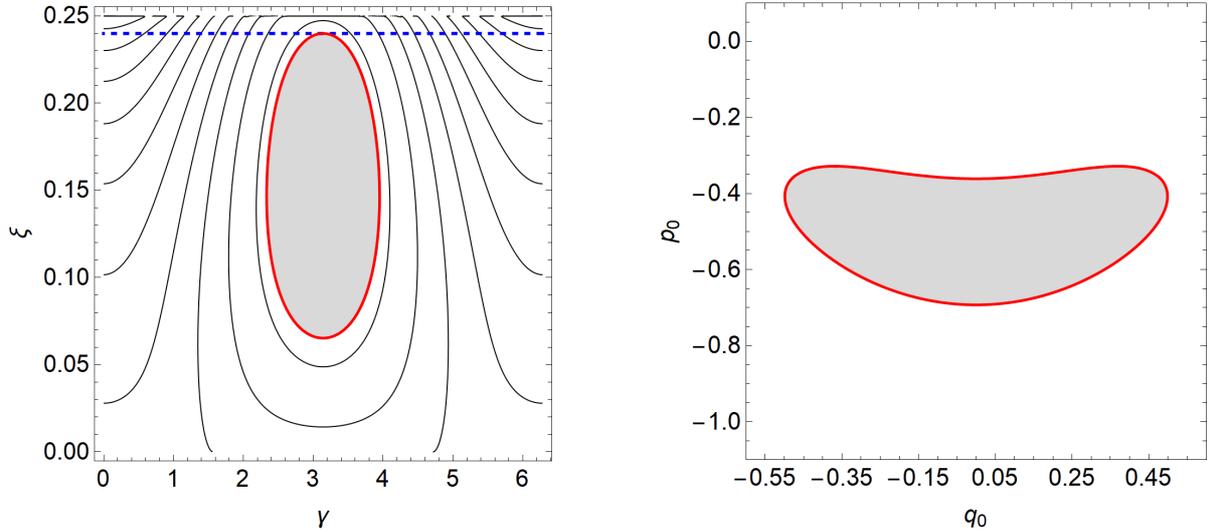

**Figure 5**: Grey region is SBMT of island kind on the $(\gamma,\xi)$ cylinder (left) and on the $(q_0,p_0)$ plane (right). Thick dashed line corresponds to truncation level $\xi_{max}=0.24$. The parameters of external forcing are $F=0.0876$; $\Omega=0.92$.

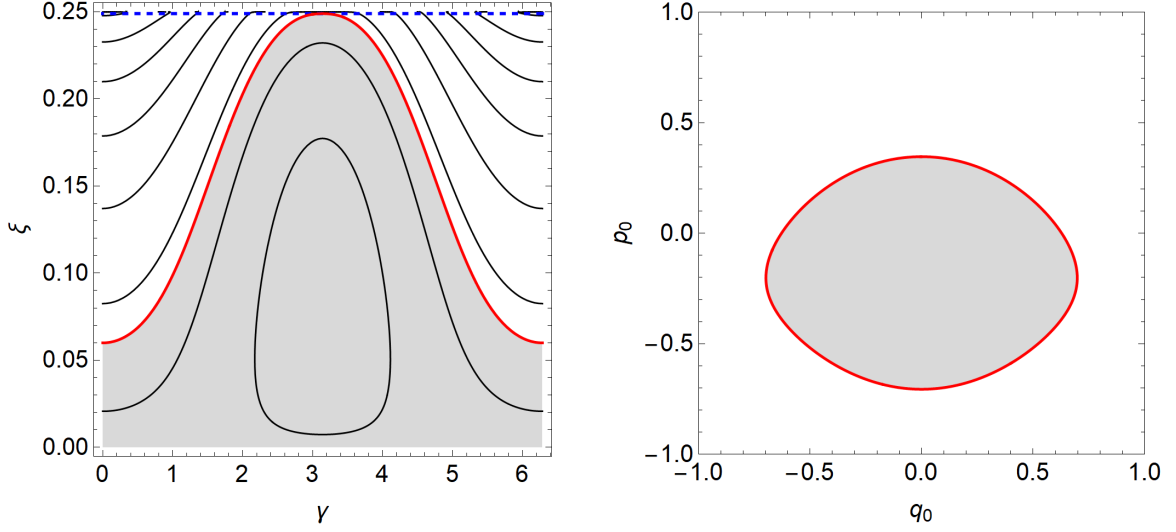

**Figure 6**: Peninsula SBMT for the parameters F=0.0876, Ω=1.04.
Thick dashed line corresponds to truncation level ξ_max=0.249.

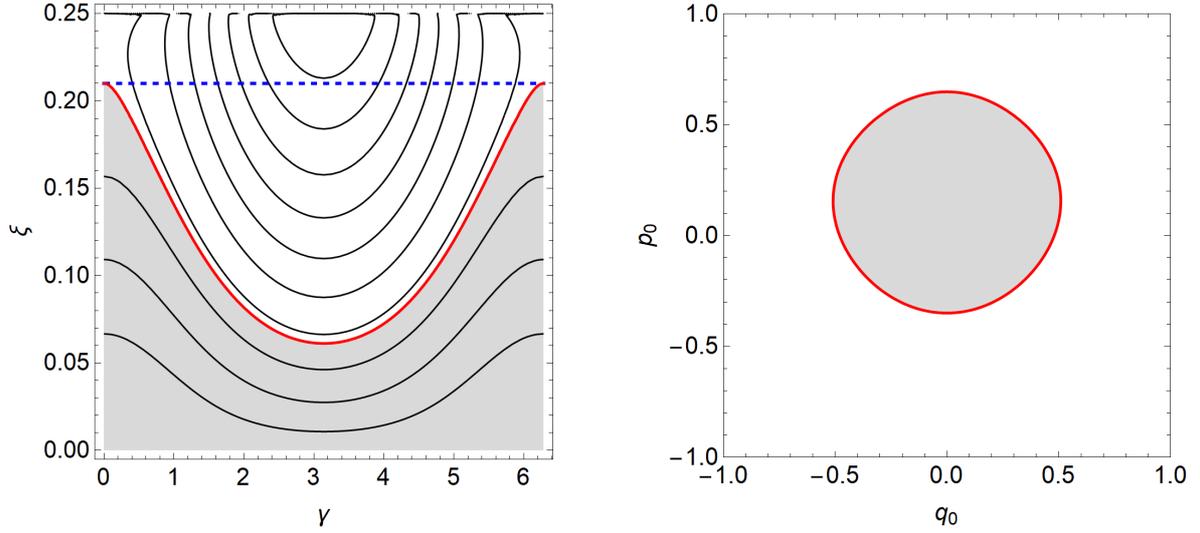

**Figure 7**: Peninsula SBMT for the parameters
F=0.0876, Ω=0.57. Thick dashed line corresponds to truncation level ξ_max=0.21.

The boundary of SBST is a phase curve which passes through the saddle point, i.e., the curve defined by the equation

$$C(\gamma, \xi) = C(\gamma^\dagger, \xi^\dagger), \tag{28}$$

where $\xi^\dagger$ and $\gamma^\dagger$ satisfy equation (25). The safe region is bounded by the curve (28) and the circle $\xi = 0$.

An example of SBST is presented in Figure 8. The case of coexisting regions is demonstrated in Figure 9. We adopt the same notation as in Figures 5 – 7.

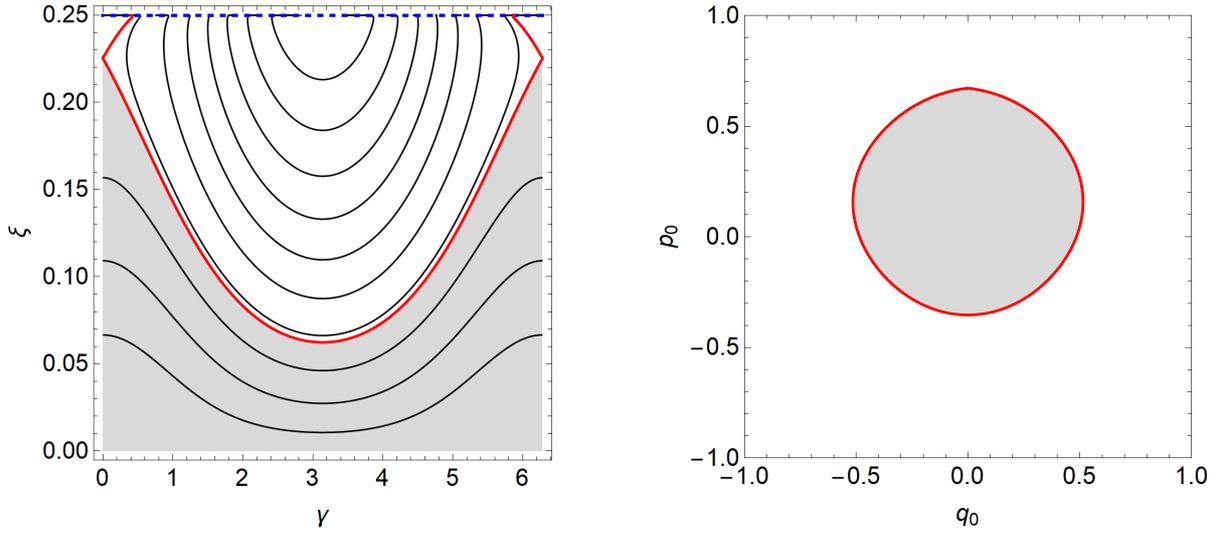

**Figure 8**: SBST for the parameters F=0.0876, Ω=0.57.
Thick dashed line corresponds to truncation level $\xi_{max}$ =0.2499.

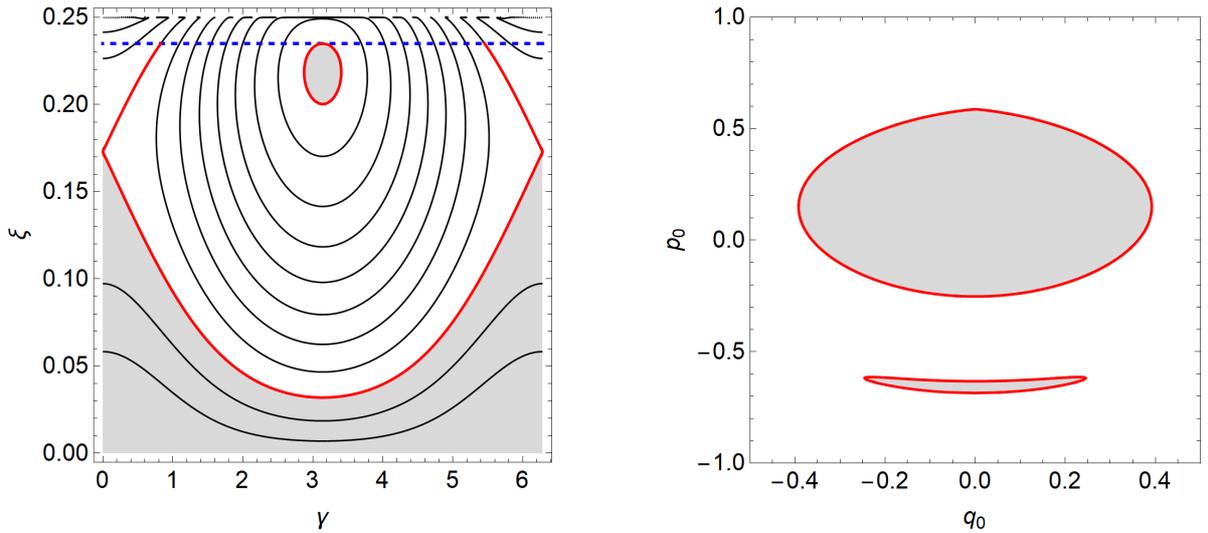

**Figure 9**: Co-existence of two types of SB with parameters F=0.0478, Ω=0.76.
Thick dashed line corresponds to truncation level $\xi_{max}$=0.235.

The amount and the location of the available critical points in equation (18) exhaust all the possibilities of escape scenarios. Therefore, for the problem in question the presented classification of SBs is complete.

## 3. Numeric verification

In order to validate theoretical predictions of the critical forcing curve $F_{cr}(\Omega)$ as well as the SBs, we compare the analytic results from Sections 1–2 with numerical simulations. The averaged system (18) is expected to approximate the original system (1) in the vicinity of the primary 1:1 resonance for relatively small excitation amplitudes. However, for the larger forces we anticipate SB to undergo erosion due to the higher order resonances, see [2], [27]–[32].

## 3.1. Critical escape force in the control plane

Figure 10 demonstrates the comparison of theoretical (solid) critical force curve with numerical simulations (dots). Every dot corresponds to a numerical estimation of $F_{cr}$ using the bisection method with accuracy of $5·10^{-5}$. Each simulation is performed for up to 3,000 *excitation cycles* (EC), i.e., a trajectory is considered non-escaping if the escape does not occur within the time interval [0,3000EC]. Three subplots of Figure 10 represent different energy truncation levels.

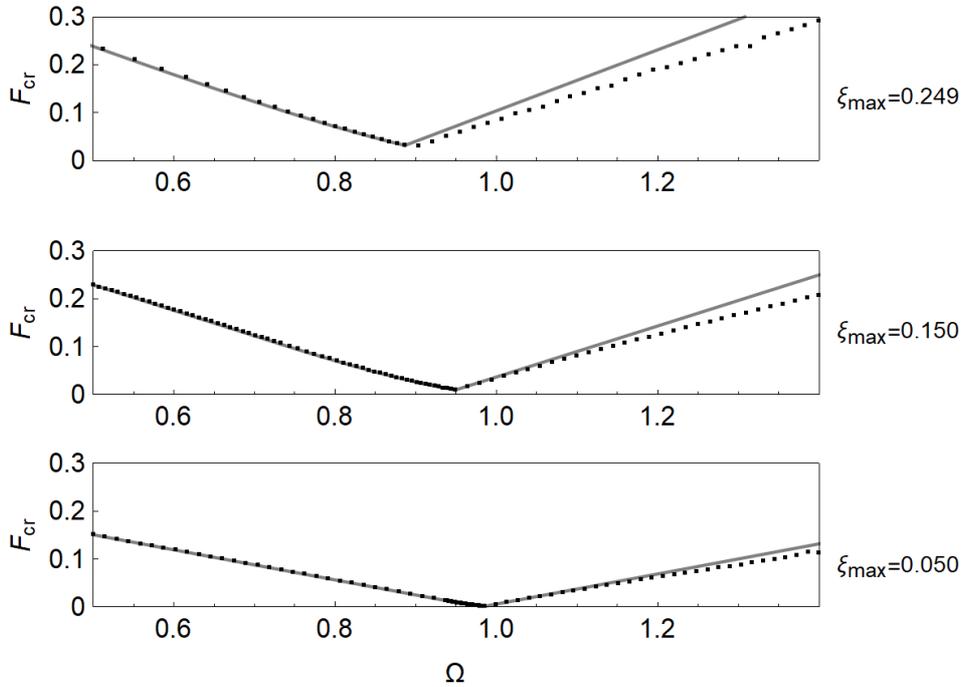

**Figure 10**: Escape threshold versus excitation frequency for different energy truncation levels $\xi_{max}$ and the zero IC. Dots represent values obtained numerally, while solid curves correspond to the theoretical prediction.

Similarly, Figures 11 and 12 show the corresponding curves $F_{cr}(\Omega)$ for nonzero ICs.

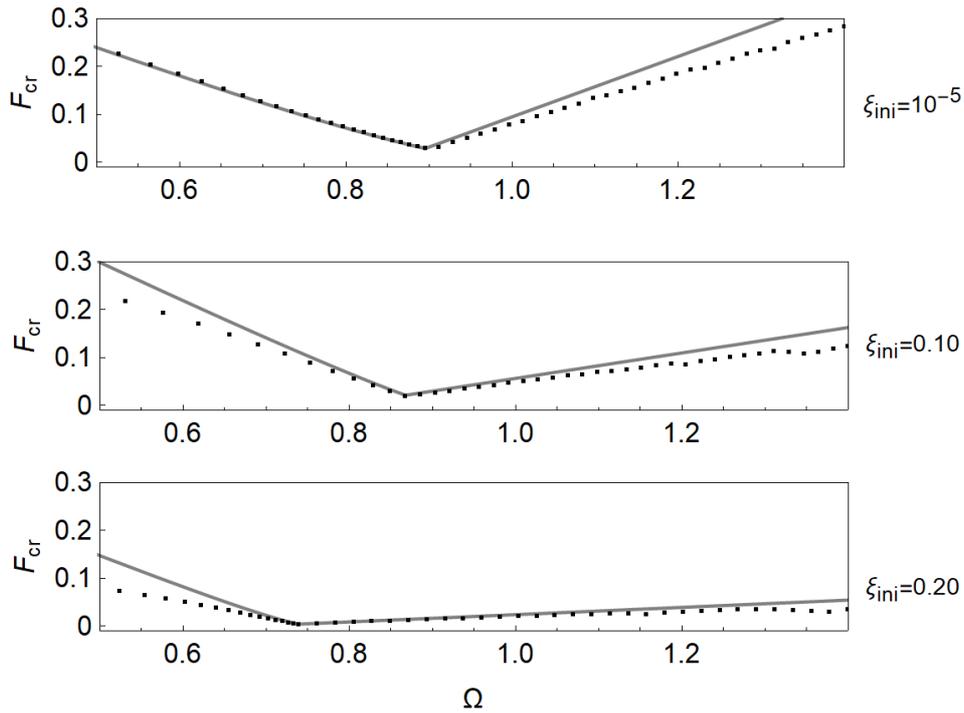

**Figure 11**: Escape threshold versus excitation frequency for various values of $\xi_{ini}$ and $\xi_{max}$ =0.242, $\gamma_{ini}$ =0.25, $\Psi=\pi$. The notation is the same as in Figure 10.

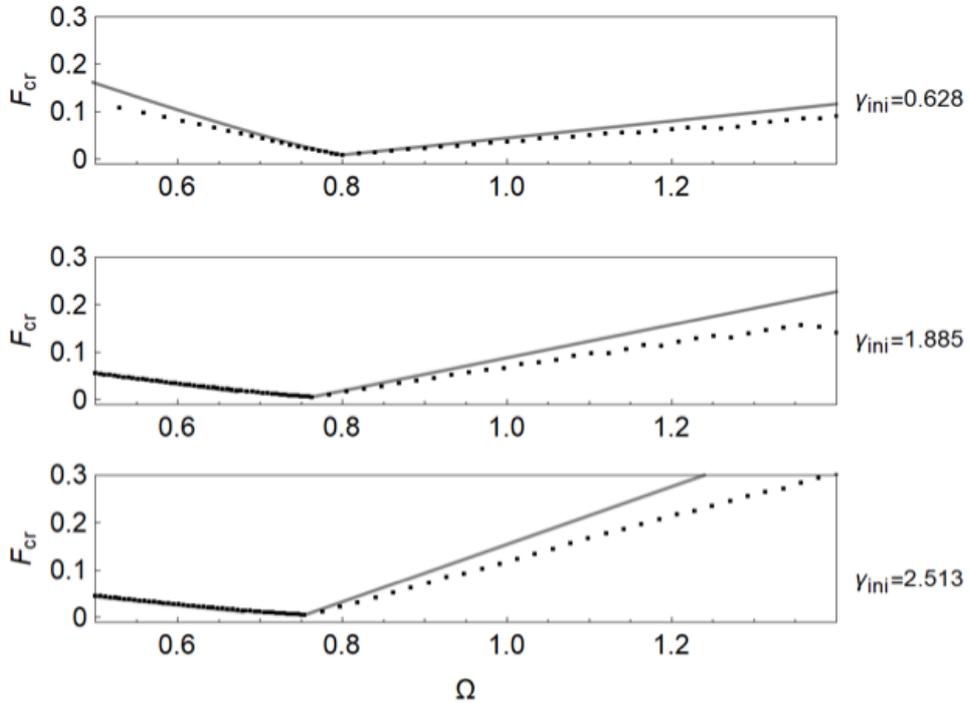

**Figure 12**: Escape threshold versus excitation frequency for various values of $\gamma_{ini}$ and $\xi_{max}$ =0.242, $\xi_{ini}$ =0.15, $\Psi=\pi$. The notation is the same as in Figure 10.

Comparison of numerical simulation results with these derived from a simplified model prediction clearly shows high fidelity approximation in the vicinity of the main resonance, the dip, and weaker prediction quality on the distant frequencies, especially with the MM scenario.

## 3.2. The SB – analytic predictions versus numerical simulations

In order to verify the approximation of SB, we compare it with the numerical simulation results. We superimpose theoretic SB boundaries on the non-escaping regions obtained with numeric time integration of system (1). The saddle point $\left(\gamma^{\dagger},\xi^{\dagger}\right)$ required for determining the boundary of an SBST is approximated using Newton-Raphson method. For all the numerical experiments the ICs are drawn from a uniform 200×200 grid.

Figure 13 and Figure 14 present the variations of the SBs with increasing energy truncation level $\xi_{max}$. Numerically obtained SBs are shaded according to the value of $\xi_{max}$. Lighter grey regions (smaller $\xi_{max}$) include all the darker ones (larger $\xi_{max}$). White colour corresponds to the escape region common to all considered values of the parameters. Each curve in Figure 13–Figure 17 represents a theoretical prediction of the SB boundary for the corresponding parameter values (see legends for the detailed description of the patterns). While the demonstrated SBs in Figure 13 have the same topology for all the cutoff energy levels, Figure 14 shows the disjoint SBs for smaller values of $\xi_{max}$ (cf. Figure 9). For the chosen parameters the agreement between theoretically obtained SBs and their numerical counterparts is very good.

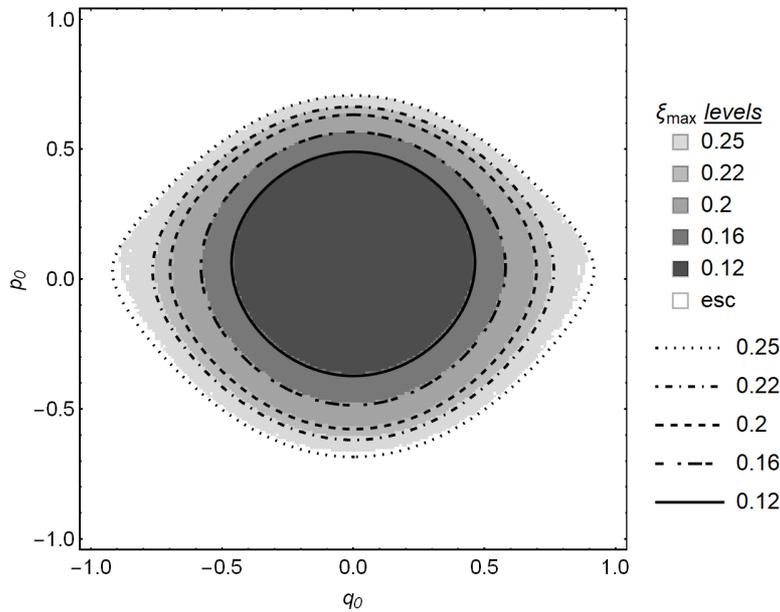

**Figure 13**: Transformation of the SB with increasing energy truncation level $\xi_{max}$. Grey regions represent the numerically obtained SBs, curves correspond to the theoretical predictions. The other parameters are F=0.015, Ω=1.04, Ψ=π.

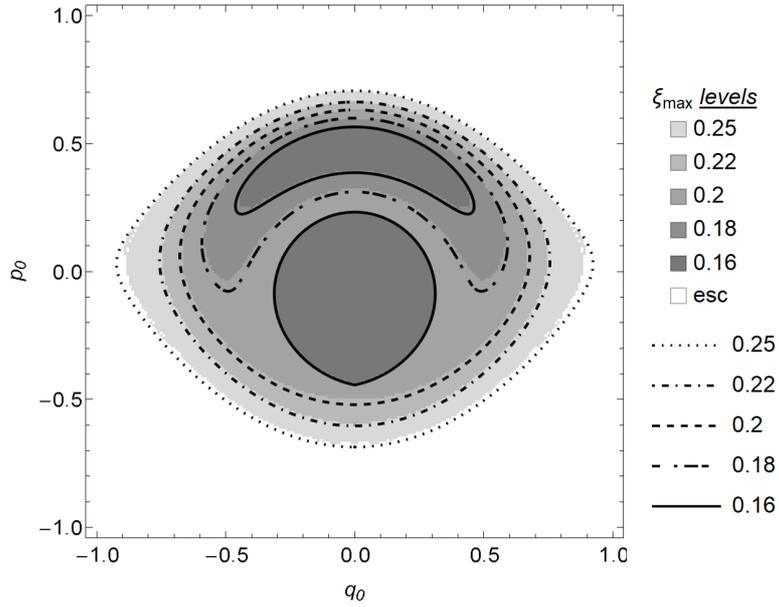

**Figure 14**: Same as Figure 13 but for F=0.01; Ω=0.9; Ψ=π;

Distortion and rotation effects on SB boundary shapes caused by the phase shift of the forcing are demonstrated in Figure 15 (cf. Figure 14).

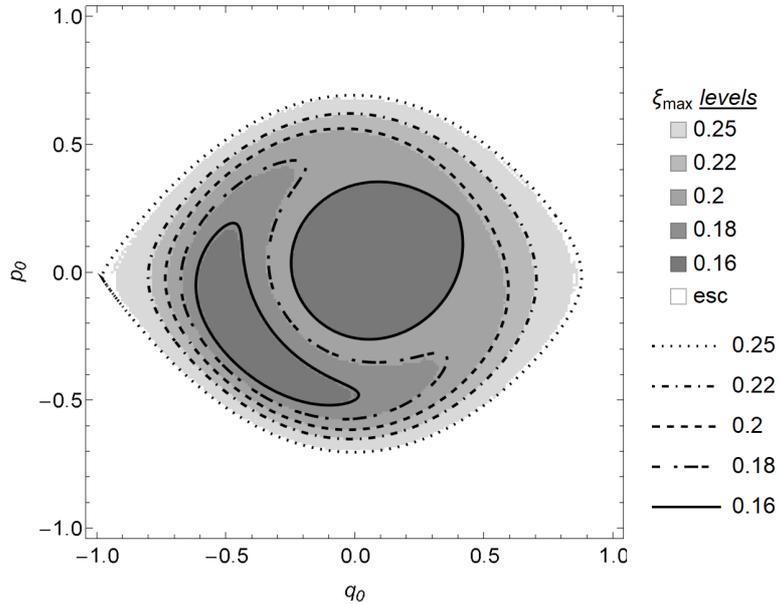

**Figure 15**: Same as Figure 14 but for Ψ=1.

For relatively small force amplitudes (Figure 16) the numeric results nicely conform to the theoretical predictions. Figure 17 presents a similar comparison for higher levels of forcing.

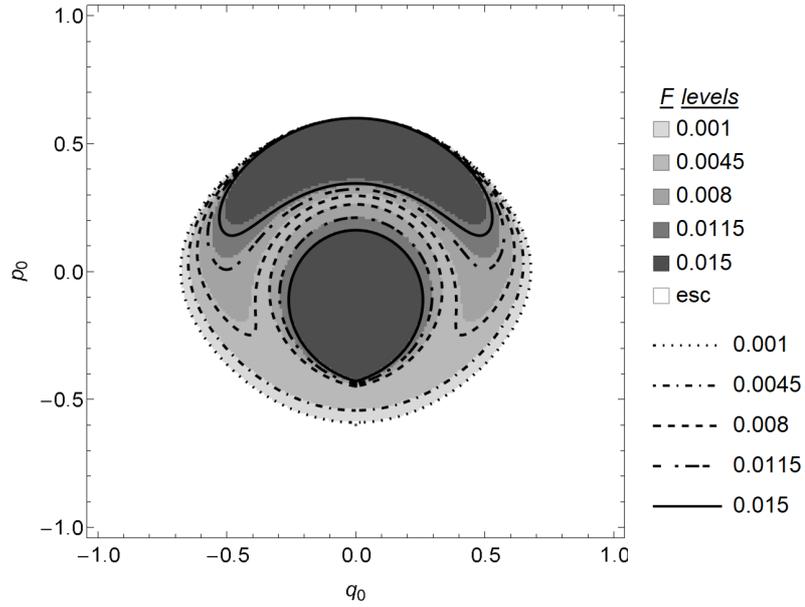

**Figure 16**: Transformation of the SB with increasing external forcing amplitude $F$. Other parameters are $\xi_{max}$=0.18, $\Omega$=0.9, $\Psi$=$\pi$.

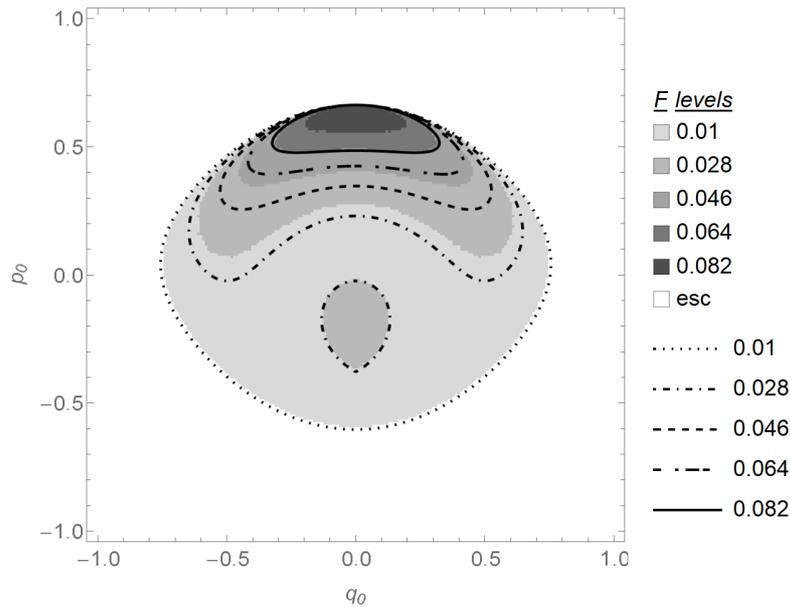

**Figure 17**: Same in Figure 16 but for $\xi_{max}$=0.22, $\Omega$=0.9, $\Psi$=$\pi$.

As one can see, for small values of the external forcing $F$ the analytically predicted boundaries perfectly match the numeric SBs, even when the SB is disjoint. However, when $F$ increases, the discrepancy between the theoretic prediction and the numeric SBs becomes more and more substantial. The reason for the reduction of the prediction quality can be attributed to the fact that the AIR method captures only the primary 1:1 resonance, while neglecting the secondary resonances which cause the erosion of SBs [27]–[32]. Nevertheless, the obtained theoretical boundaries can serve as a reasonable initial approximation of the SB for subsequent possible numeric refining.

In order to further analyze the structure of SB for the larger values of external forcing, it is instructive to consider the stroboscopic (period) map

$$(q(t), p(t)) \mapsto (q(t+T), p(t+T)), \qquad (29)$$

with $T = 2\pi/\Omega$. Figure 18 presents the phase portraits of 100 non-escaping trajectories of the map (29) starting from random ICs. We observe one big island of invariant tori centred around a fixed point, many smaller islands of invariant tori centred around various periodic solutions, and a stochastic layer in between. The structure of the phase portrait of the map (29) suggests the presence of the higher order resonances contributing to the erosion of the SB. For example, invariant tori around a 6-periodic solution presented on the figure, indicate the existence of 1:6 resonance. Dashed curve represents the theoretical prediction of the SB boundary, and although it fails to capture the SB exactly, it still serves as a reasonable first guess for locating the SB.

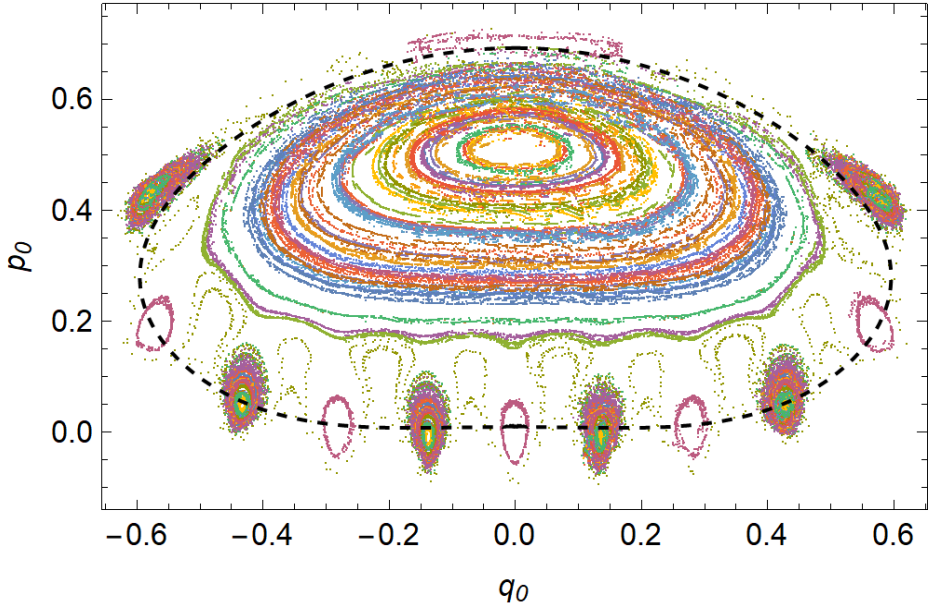

**Figure 18**: Phase portrait of 100 non-escaping trajectories of the period map (28) for the following set of parameters: F=0.0876, Ψ=π, Ω=0.99, $\xi_{max}$=0.24. Each trajectory runs for 3000 iterations. Dashed curve corresponds to the theoretically predicted boundary of the SB.

## 4. Discussion and conclusions

The results obtained above indicate that the approximation of isolated resonance enables analytic prediction of the safe boundaries in the problem of forced escape. For relatively small forcing amplitudes, these predictions are very accurate. The reason is that the forcing term is in fact a perturbation of integrable basic Hamiltonian describing the conservative motion in the well. This approach reveals some new possible properties of the SB, e.g., complicated, or even disjoint geometric shape

For higher forcing amplitudes, the quality of the prediction deteriorates due to intrinsic inaccuracies of the AIR and, more significantly, due to secondary resonances that are not taken into account in the AIR and cause the SB erosion. Still, the analytic prediction can serve as a viable initial guess for the SB with further numeric tightening.

The analysis presented above is substantially based on the ability to perform the action-angle transformation in the benchmark model with quartic potential. Application of the developed method to more generic models, where this transformation is not possible, may be based on appropriate approximation techniques, and requires further exploration.

**Acknowledgement.**

The authors are very grateful to Israel Science Foundation (grant 1696/17) for financial support.

# Appendix: Comparison of the escape criteria

As it was mentioned before, the two escape criteria — maximum energy and maximum displacement — are not equivalent. It is obvious that the energy criterion immediately follows from the displacement criterion, and therefore, the SB based on the *energy criterion* (ESB) is always included in the corresponding SB defined through the maximum *displacement criterion* (DSB). However, the other way is not necessarily true.

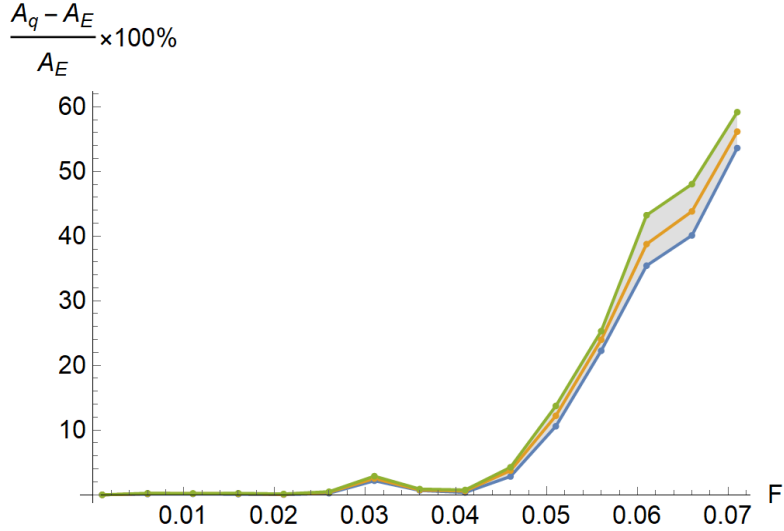

**Figure 19**: Comparison of SB areas for different values of the external force amplitude F. Blue, orange and green markers correspond to the minimum, the mean, and the maximum of 5 experiments. Other parameters are the following: Ω=0.9, Ψ=0, $\xi_{max}$=0.25.

In order to compare two criteria, the following numeric test is performed. For each value of the forcing amplitude *F*, we pick 10,000 random ICs from the square $[-1,1] \times [-1,1]$. Then, we perform numerical time integration of system (1) starting from the chosen set selecting only the non-escaped trajectories according to the ESB. Next, we run the simulations again for the rest of ICs but with the displacement criterion to obtain the set belonging to the DSB but not to the ESB. Let $A_q$ and $A_E$ denote the number of points that belong to the DSB and the ESB, respectively. The difference between $A_q$ and $A_E$ relative to $A_E$ exhibits the discrepancy between the two criteria.

The experiment is repeated 5 times for the values of the external forcing amplitude *F* ranging from 0.001 to 0.071 with a step $\Delta F = 0.005$. The results are presented in Figure 19. As one can see, for the small values of the external forcing amplitude *F*, SBs for both criteria are almost identical. However, with increasing *F* the difference becomes substantial, see Figure 19.

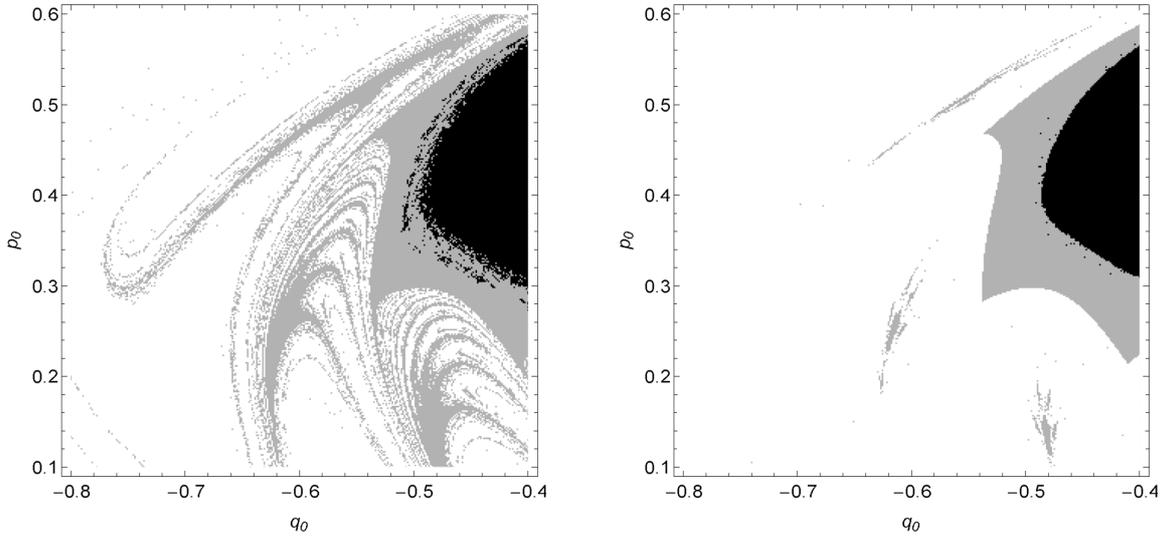

**Figure 20**: Zoomed SB Escape Chart simulation for
F=0.0876, Ω=0.95, Ψ=π, ξ$_{max}$=0.25, the resolution is 300×300.
Left: t$_{eval}$=500EC. Right: t$_{eval}$=11,000EC, DSB (grey) and by ESB (black).

Furthermore, our findings reveal that the choice of the escape criterion also strongly affects the shape of the SBs. The maximum energy escape criterion yields SBs with smoother boundary. Figure 20 presents a comparison of SBs in the ($q_0$,$p_0$) plane for two superimposed escape criteria for a short evaluation time $t_{eval}$ (left) and a long $t_{eval}$ (right). Colours grey and black correspond to the displacement and the energy criteria, respectively.

For a short evaluation time $t_{eval}$ the boundaries of SBs for both criteria look irregular representing the transient behaviour. For a long time $t_{eval}$ after the transient dynamics fades out, the SB boundary becomes clearer and smoother. The full view of SB for the same values of parameters as in Figure 20 is shown on the Figure 21.

Figure 22 presents the change of the SB area with time. As one can see the area undergoes a rapid shrinking for about 2000EC, followed by almost no change after. Each point on the graph is a natural logarithm of the SB area computed numerically at a given evaluation time.

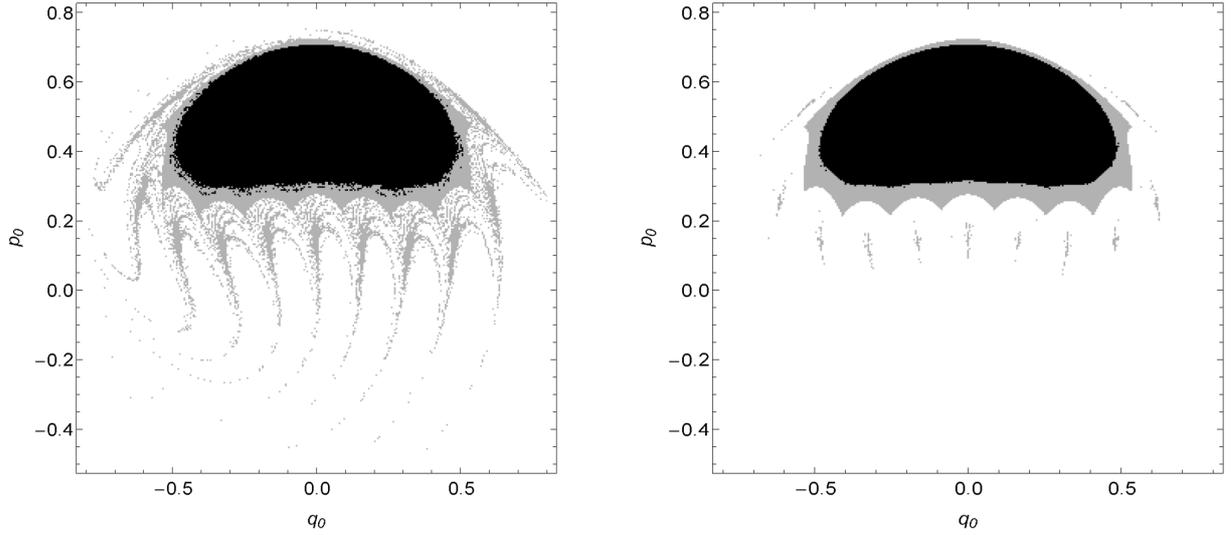

**Figure 21**: Entire SB Escape Chart simulation for
F=0.0876, Ω=0.95, Ψ=π, $\xi_{max}$=0.25, the resolution is 300×300.
Left: $t_{eval}$=500EC. Right: $t_{eval}$=11,000EC, DSB (grey) and by ESB (black).

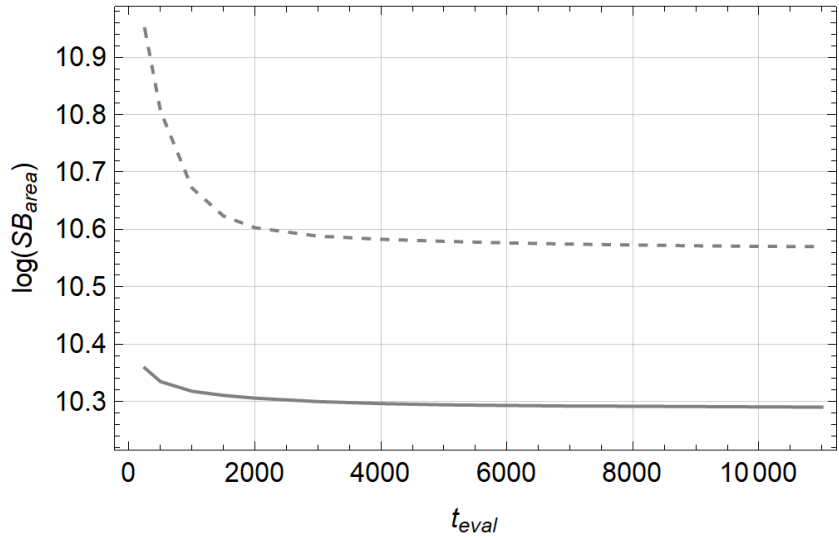

**Figure 22**: Log of SB area in pixels fading with $t_{eval}$ for DSB (dashed) and ESB (solid)
F=0.0876, Ω=0.95, Ψ=π, $\xi_{max}$=0.25, the resolution is 300×300.

26